\documentclass[12pt,journal,final,finalsubmission,onecolumn]{IEEEtran}
\usepackage{graphicx}
\usepackage{subfigure}
\usepackage{amsmath,amssymb,color}
\usepackage{cite}

\title{Global alignment of protein-protein interaction networks by graph matching methods}
\author{Mikhail Zaslavskiy,
        Francis Bach,
        and~Jean-Philippe~Vert
\thanks{Mikhail Zaslavskiy is with the Centre for Computational Biology and the Centre for Mathematical Morphology; Jean-Philippe Vert is with the Centre for Computational Biology, Mines ParisTech, 35 rue Saint-Honor\'e, 77305 Fontainebleau, France. They are also with the Institut Curie, Paris, F-75248 France, and with INSERM U900, Paris, F-75248 France.}
\thanks{E-mail: mikhail.zaslavskiy@ensmp.fr, jean-philippe.vert@ensmp.fr}
\thanks{Francis Bach is with the INRIA-Willow Project Team, Laboratoire d'Informatique de l'Ecole Normale Sup\'erieure 
(CNRS/ENS/INRIA UMR 8548), 45 rue d'Ulm, 75230 Paris, France}
\thanks{E-mail: francis.bach@mines.org}}

\newcommand {\todo}[1]{{\bf TODO: #1}}

\newcommand{\tr}{ {\rm tr }}

\newcommand{\Pcal}{\mathcal{P}}

\newcommand{\Acal}{\mathcal{A}}
\newcommand{\Tcal}{\mathcal{T}}
\newcommand{\RR}{\mathbb{R}}
\newcommand{\OMIT}[1]{}

\begin{document}

\maketitle

\begin{abstract}
\noindent \textbf{Motivation:}\\
Aligning protein-protein interaction (PPI) networks of different species has drawn a considerable interest recently. This problem is important to investigate evolutionary conserved pathways or protein complexes across species, and to help in the identification of functional orthologs through the detection of conserved interactions. It is however a difficult combinatorial problem, for which only heuristic methods have been proposed so far.\\
\noindent \textbf{Results:}\\
We reformulate the PPI alignment as a graph matching problem, and investigate how state-of-the-art graph matching algorithms can be used for that purpose. We differentiate between two alignment problems, depending on whether strict constraints on protein matches are given, based on sequence similarity, or whether the goal is instead to find an optimal compromise between sequence similarity and interaction conservation in the alignment. We propose new methods for both cases, and assess their performance on the alignment of the yeast and fly PPI networks. The new methods consistently outperform state-of-the-art algorithms, retrieving in particular $78\%$ more conserved interactions than IsoRank for a given level of sequence similarity.\\
\noindent \textbf{Availability:}\qquad {\tt http://cbio.ensmp.fr/proj/graphm\_ppi/}\\
Additional data and codes are  available upon request.\\
\noindent \textbf{Contact:} {\tt jean-philippe.vert@mines-paristech.fr}
\end{abstract}

\section{Introduction}

PPIs play a central role in most biological processes. Recent years have witnessed impressive progresses towards the elucidation of large-scale PPI networks in various organisms, thanks in particular to the development of high-throughput experimental techniques such as yeast two-hybrid \cite{Fields1989novel} or coimmunoprecipitation followed by mass-spectrometry \cite{Aebersold2003Mass}. As the amount of PPI network data increases, computational methods to analyze and compare them are also being developed at a fast pace. In particular, comparative PPI network analysis across species has already provided insightful views of similarities and differences between species at the systemic level \cite{Sharan2005Conserved,Suthram2005Plasmodium} and helped in the identification of functional orthologs \cite{Bandyopadhyay2006Systematic}.

Comparing PPI networks usually involves some form of \emph{network alignment}, i.e.,  the identification of pairs of homologous proteins from two different organisms, such that PPIs are conserved between matched pairs. The rationale behind this notion is that a protein and its functional orthologs are likely to interact with proteins in their respective network that are themselves functional orthologs. Hence, while direct sequence homology alone is often not sufficient to identify functional orthologs within paralogous families \cite{Sjoelander2004Phylogenomic}, the use of PPI information can help in the disambiguation of functional orthologs within clusters of homologous sequences, such as those produced by the Inparanoid algorithm \cite{Remm2001Automatic}. This approach has been investigated in particular by \cite{Bandyopadhyay2006Systematic}. Conversely, network alignment can also be a valuable approach to validate PPI conserved across multiple species and detect evolutionary conserved pathways or protein complexes \cite{Sharan2005Conserved,Kelley2003Conserved}.

Several methods have been proposed to perform \emph{local} network alignment (LNA) of PPI networks, i.e., to identify subsets of matching pairs of proteins with conserved subgraphs of interactions. These methods include PathBLAST \cite{Kelley2003Conserved,Kelley2004PathBLAST} and NetworkBLAST \cite{Sharan2005Conserved}, which adapt the ideas of the BLAST algorithm to the search for local alignments between graphs, the method of \cite{Koyutuerk2006Pairwise}, inspired by biological models of deletion and duplication, Graemlin \cite{Flannick2006Graemlin}, which uses networks of modules to infer the alignment, or the Bayesian approach of \cite{Berg2006Cross-species}. Less attention has been paid to the problem of \emph{global} network alignment (GNA), i.e., the search for a global correspondence between most or all vertices of two networks which again matches similar proteins and leads to conserved interactions. Notable exceptions include the Markov random field (MRF) based method of \cite{Bandyopadhyay2006Systematic} and the IsoRank algorithm \cite{Singh2008Global} which formulates the problem as an eigenvalue problem.

While LNA procedures can detect multiple, unrelated matched regions between networks, and can in particular match a given protein of a network to several proteins of the other network in different local matchings, GNA seeks the best consistent matching across all nodes simultaneously. This can be a desirable property for many applications, such as functional ortholog identification. On the other hand, from a computational point of view, GNA is arguably more difficult than LNA since it must find a solution among all possible global matchings. In fact, as we explain below, it is natural to reformulate GNA as weighted graph matching problem, a problem for which no polynomial-time algorithm is known. Solving the general GNA problem therefore must involve some sort of approximate or heuristic method, such as IsoRank.

Following this line of thought, we propose here to formulate explicitly GNA as a graph matching problem, and investigate the use of modern state-of-the-art exact and approximate methods to solve it.  While no exact solution of the graph matching optimization problem can be found in general, we show that in certain cases, if ``enough constraints'' are put on the possible protein associations, and if the PPI networks are ``not too dense'' (these notions being rigorously defined in Section \ref{sec:optimal}), then an exact solution can be found efficiently by a new message-passing algorithm. Interestingly, this case arises in particular in the functional ortholog detection problem between yeast and fly investigated by \cite{Bandyopadhyay2006Systematic}, where matching pairs are constrained to belong to clusters of proteins produced by the Inparanoid algorithm and the PPI networks of both species are not too dense. On these data, we are therefore able to find a matching which conserves more interactions than the solutions found by MRF \cite{Bandyopadhyay2006Systematic} as well as a version of  IsoRank adapted to this situation \cite{Singh2008Global}, and we are in fact certain that our solution is optimal in the sense that it produces the largest possible number of conserved interactions. Interestingly, the resulting alignment retrieves $13\%$ more HomoloGene pairs than the alignments of MRF and $5\%$ more than that of IsoRank, suggesting that maximizing the number of conserved interactions indeed improves functional orthology disambiguation. When the GNA is more complex, e.g., matched pairs are not limited to belong to the same Inparanoid clusters, or the PPI networks have more edges, then our message-passing algorithm can not be used and the optimal matching can not be found in reasonable time anymore. In that case we propose to use a recent state-of-the-art approximate methods for graph matching \cite{Zaslavskiy2008path}, which tracks a path of solutions for a family of relaxed problems, as well as a new, faster and more direct gradient-based method, which bears similarities with the IsoRank method. Like IsoRank, these methods have a free parameter to balance the trade-off between matching similar proteins, on the one hand, and producing an alignment with many conserved interactions, on the other hand. We test them on the global unconstrained alignment of the fly and yeast networks, and show that for a given level of mean sequence similarity between matched proteins, our new method retrieves $ 78\% $ more conserved interactions than IsoRank.

\section{Constrained and balanced GNA problems}
\label{sec:problem}
In this section we set the notations and formalize two variants of the GNA problems.
We represent a PPI network describing the interactions among $N$ proteins of an organism as an undirected simple graph $G=(V_G,E_G)$, where $V_G=\left(v_1,\ldots,v_N\right)$ is a finite set of $N$ vertices representing the $N$ proteins, and $E_G \subset V_G \times V_G$ is the set of edges representing the pairs of interacting proteins. Each such graph (or network)  can equivalently be represented by a symmetric $N\times N$ adjacency matrix $A_G$  where $[A_G]_{ij}=[A_G]_{ji}=1$ if protein $v_i$ interacts with protein $v_j$ and $0$ otherwise.

Given two graphs $G$ and $H$ representing the PPI networks of two species, the GNA problem is, roughly speaking, to find a correspondence between the vertices of $G$ and the vertices of $H$ which matches similar proteins and enforces as much as possible the conservation of interactions between matched pairs in the two graphs. To formalize this, let us assume that $G$ and $H$ have the same number $N$ of vertices, and that we are looking for a bijection between the vertices of $G$ and the vertices of $H$. Although this may sound at first sight a strong assumption, given that PPI networks usually do not have the same size, and that we may not want to match all proteins of each network, both limitations can be addressed by adding dummy nodes (with no connection) to each graph in order to ensure that they finally have the same size. In a complete matching of such graphs with dummy nodes, matching a protein to a dummy node simply means that in the GNA the protein is not matched. $G$ and $H$ being assumed to have the same number of vertices, a matching of their vertices is now simply a permutation $\pi$ of $\{1,\ldots,N\}$ which associates the $i$-th vertex of $H$ with the $\pi(i)$-th vertex of $G$. Equivalently, the permutation $\pi$ can be represented by a $N\times N$ permutation matrix $P$, i.e., a binary matrix whose $(i,j)$-th entry is equal to $1$ if and only if $\pi(i)=j$ (that is, when the $i$-th vertex of $H$ is matched to the $j$-th vertex of $G$). We denote by $ \Pcal =\{P\in \{0,1\}^{N\times N}:P1_N=1_N, P^T1_N=1_N\}$ the set of permutation matrices, where $1_N$ is the $N$-dimensional vectors whose entries are all equal to $1$.

The number of interactions conserved by a permutation $\pi$ is the number of pairs $(i,j)$ which are connected in $H$, and such that their corresponding vertices $\pi(i)$ and $\pi(j)$ are also connected in $G$. Let us denote by $J(P)$ the number of such interactions conserved by the permutation encoded in the permutation matrix $P$. In order to express $J(P)$, we can observe that if we apply the permutation encoded by $P$ to the vertices of $H$, we obtain a new graph isomorphic to $H$ which we denote by $P(H)$. It is easy to see that the adjacency matrix of the permuted graph, $A_{P(H)}$, is simply obtained from $A_H$ by the equality $A_{P(H)}=PA_HP^T$ \cite{Umeyama1988eigendecomposition}. As a result, $J(P)$ is simply obtained as half the number of entries which are simultaneously equal to $1$ in both binary matrices $A_G$ and $PA_HP^T$ (each conserved interaction results in two identical entries, by symmetry of the adjacency matrices). Hence we obtain the following expression for $J(P)$:
\begin{equation}
\label{eq:JP}
J(P) = \frac{1}{2}\sum_{i,j=1}^N [A_G]_{ij} [PA_HP^T]_{ij} = \frac{1}{2}\tr(A_G^TPA_HP^T)\,.
\end{equation}

Besides the number of conserved interactions, a good GNA should match proteins with similar sequences. We consider here two possible formulations of this objective.
\begin{itemize}
\item \emph{Constrained GNA}. Here we assume that a pre-processing of the protein sequences has produced a set of candidate matched pairs $\Acal \subset V_H \times V_G$, and we simply wish to disambiguate the matching using PPI information, if some proteins have several candidate matchings. This is for example the formulation proposed by \cite{Bandyopadhyay2006Systematic}, where a first clustering of all proteins sequences is performed to define a collection of protein clusters with the Inparanoid algorithm, and the pairs matched between the yeast and fly proteome are constrained to belong to the same cluster. Such constraints can be directly encoded as constraints over the permutation matrix $P$, by imposing $P_{ij}=0$ if the $i$-th vertex of the first graph and the $j$-th vertex of the second graph are not allowed to match. We are then looking for a solution in the set of matrices $\Pcal_\Acal = \left\{P\in\Pcal\,:\,\forall (i,j)\in [1,N]^2 \backslash \Acal ,P_{ij}=0 \right\}$, and it is then natural to look for the permutation compatible with the constraints with the largest number of conserved interactions, i.e., to solve:
\begin{equation}\label{eq:constrainedGNA}
\max_{P\in\Pcal_\Acal}J(P)\,.
\end{equation}
\item \emph{Balanced GNA}. A interesting property of constrained GNA is that, by reducing the search space to $\Pcal_\Acal$, it can result in a tractable optimization problem (as shown for example in Section \ref{sec:optimal}). On the other hand, in some cases one may want to accept matching between less similar vertices if it leads to an important increase in the number of conserved interactions. In other words, one would like to be able to automatically \emph{balance} the matching of similar vertices with the conservation of interactions, as advocated by \cite{Singh2008Global} and implemented by IsoRank. This can be formalized by assuming that a $N\times N$ matrix of similarities between vertices $C$ is given (e.g., derived from pairwise sequence similarity scores), and by trying to maximize the total similarity between matched pair. $C_{ij}$ denoting the similarity between the $i$-th vertex of $G$ and the $j$-th vertex of $H$, the total similarity between pairs matched by a permutation matrix $\pi$ is simply
\begin{equation}\label{eq:SP}
S(P) = \sum_{i=1}^N C_{\pi(i),i} = \tr\left(PC\right)\,.
\end{equation}
In order to find a balance between matching similar pairs (large $S(P)$) and having many conserved interactions (large $J(P)$), we propose to consider the following optimization problem:
\begin{equation}\label{eq:balancedGNA}
\max_{P\in\Pcal} \lambda J(P) + (1-\lambda)S(P)\,,
\end{equation}
where $\lambda\in[0,1]$ controls the trade-off between both objectives. $\lambda=1$ corresponds to the maximization of $J(P)$ only, i.e., to find a good topological matching of the graphs independently of the similarity between matched pairs, while $\lambda=0$ amounts to focus only on the similarity between proteins and finding a matching which maximized the mean sequence similarity, without using PPI information.
\end{itemize}
When $\lambda>0$, the balanced GNA problem (\ref{eq:balancedGNA}) is equivalent to a general graph matching problem, discussed in Section \ref{sec:approximate}, which is known to be computationally intractable in general. The constrained GNA (\ref{eq:constrainedGNA}) can be seen as a particular case of the balanced GNA, by taking the similarity function equal to $0$ between two vertices allowed to match and $-\infty$ for two vertices not allowed to match. Indeed, in that case (\ref{eq:balancedGNA}) is equivalent to minimizing $J(P)$ over the set of matrices $P$ for which $S(P)$ is finite, that is exactly the set $\Pcal_\Acal$ of (\ref{eq:constrainedGNA}). While indeed general graph matching methods to solve (\ref{eq:balancedGNA}) can be applied to solve (\ref{eq:constrainedGNA}), we show in the next Section that in some cases there exists a simple polynomial-time algorithm to solve (\ref{eq:constrainedGNA}) directly even for large non-sparse graphs.

\section{Methods}
In this section we present methods to solve both the constrained GNA problem (\ref{eq:constrainedGNA}) and the balanced GNA problem (\ref{eq:balancedGNA}). Since any algorithm to solve the balanced GNA problem can also solve the constrained GNA, as explained in the previous section, we start by describing methods to solve the balanced GNA problem.

\subsection{Algorithms for the balanced GNA problem}
\label{sec:approximate}
\OMIT{In many cases we may not be able to use the exact optimization method presented in Section \ref{sec:optimal}, either because we consider the balanced GNA problem (\ref{eq:balancedGNA}), or because we consider the constrained GNA problem (\ref{eq:constrainedGNA}) but the graph of clusters investigated in Section \ref{sec:optimal} has loops and has a large tree-width.}
The balanced GNA problem (\ref{eq:balancedGNA}) is a general graph matching problem, which is known to be a difficult combinatorial problem. While some methods based on incomplete enumeration may be applied to search for an exact optimal solution in the case of small or sparse graphs, only approximate algorithms that usually find non-optimal solutions but are more scalable can be used for large non-sparse graph matching. Many such approximate algorithms have been proposed, see e.g., the review of \cite{Conte2004Thirty}. They include in particular spectral methods \cite{Umeyama1988eigendecomposition,Caelli2004eigenspace,Singh2008Global}, or methods based on a relaxation of the optimization problem (\ref{eq:balancedGNA}) \cite{Almohamad1993linear,Gold1996graduated}. They differ mainly on their scalability, and on the accuracy of the solution found. For example, a comparison of several such methods was carried out recently \cite{Zaslavskiy2008patha,Zaslavskiy2008path}. 

Based on these observation, we propose here to use state-of-the-art graph matching methods to balanced GNA for PPI networks. In particular we focus on the PATH algorithm \cite{Zaslavskiy2008path}, which was shown to provide state-of-the-art performance in various graph matching benchmark. We also propose a new and simpler gradient ascent method, similar in spirit to the Graduated Assignment (GA) algorithm \cite{Gold1996graduated}. As a benchmark, we consider the IsoRank method, which can be thought of as a particular spectral method for graph alignment, and which is currently the method of choice for balanced GNA of PPI networks. We now briefly describe these methods.
\begin{itemize}
\item \emph{PATH method.} The PATH algorithm is based on two relaxations of (\ref{eq:balancedGNA}), one concave and one convex, over the set of doubly stochastic matrices \cite{Zaslavskiy2008path}. The method starts by solving the convex relaxation, and then iteratively solves a linear combination of the convex and concave relaxations by gradually increasing the weight of the concave relaxation and following the path of solutions thus created. It finishes when the a solution reaches a corner of the set of doubly stochastic matrices, i.e., when the solution is a permutation matrix in $\Pcal$. On several benchmarks, the PATH method was shown to be state-of-the-art in accuracy, and can easily process graphs with a few thousands vertices in a few hours on a modern desktop computer.
\item \emph{GA method.} We propose a new, simple gradient method based on a relaxation of (\ref{eq:balancedGNA}) over the set of doubly stochastic matrices. Although the function to be maximized is not concave (because of the term $J(P)$), we simply start from an initial solution and iteratively choose a new permutation matrix in the direction of the gradient of the objective function. This approach may be relevant if we can start from a ``good'' initial solution, i.e., if we solve a constrained GNA (\ref{eq:constrainedGNA}) where the constraints are strong enough. The gradient of $S(P)$ in (\ref{eq:SP}) is equal to $S$, the gradient of $J(P)$ in (\ref{eq:JP}) at a matrix $P_n$ is equal to $A_G^TP_nA_H$. Hence we propose to iteratively update the permutation matrix following the rule $P_{n+1}\leftarrow \arg\max_{P \in \Pcal} \tr \left( [\lambda A_G^TP_nA_H + (1-\lambda)C ]P\right)$, which can be found efficiently by the Hungarian algorithm \cite{Kuhn1955Hungarian}.
\item \emph{IsoRank method.} 
The idea of the IsoRank algorithm is to use the following recursive formula \cite{Singh2008Global}
\begin{equation}
R(i,j)=\sum_{v\in N(i)}\sum_{u\in N(j)}\frac{1}{|N(u)||N(v)|}R(u,v),~i\in V_G,~j\in V_H,
\label{eq:R}
\end{equation}
where $N(i)$ denotes the set of neighbors of $i$, $V_G$ denotes the set of vertices of graph $G$ and element $R(i,j)$ represents the similarity between vertex $i$ of graph $G$ and vertex $j$ of graph $H$. In the case of PPI networks it represents the ``likelihood'' that proteins $i$ and $j$ are functional orthologs. The recursive formula says that the more  $i$ and $j$ have similar neighbors, the greater is the similarity measure between $i$ and $j$. To estimate $R$, \cite{Singh2008Global} propose to use the power method to iteratively update $R$ according to: 
\begin{equation}
R\leftarrow AR/||AR||\,,
\end{equation}
where $A$ is the $N^2\times N^2$ matrix defined as:
$$
A(i,j)(u,v)=\frac{1}{|N(u)||N(v)|}\,.
$$
To take into account the information on protein sequence similarities encoded by matrix $C$, the following modification of (\ref{eq:R}) is used 
\begin{equation}
R=\lambda A R+(1-\lambda)C,
\label{eq:R2}
\end{equation}
where $\lambda$ has the same interpretation as in (\ref{eq:balancedGNA}).
\end{itemize}

\subsection{Algorithms for the constrained GNA problem}
\label{sec:optimal}
As explained in Section \ref{sec:problem}, all methods for solving the balanced GNA problem (\ref{eq:balancedGNA}) can also be used to solve the constrained GNA problem (\ref{eq:constrainedGNA}), by using a particular similarity function to enforce the constraints. Hence a first series of methods to solve (\ref{eq:constrainedGNA}) are the constrained version of IsoRank, GA and PATH, described in the previous section. In addition to these three methods, we consider two additional approaches specifically dedicated to the constrained GNA problem: the Markov random field (MRF) method of \cite{Bandyopadhyay2006Systematic}, and a new method based on message passing (MP) which we propose to find the global optimum of (\ref{eq:constrainedGNA}) when the graphs are not too dense.
\begin{itemize}
\item \emph{MRF method.}
To solve ambiguous assignments in Inparanoid clusters with more than two proteins, \cite{Bandyopadhyay2006Systematic} propose to use the information on protein interactions, by choosing the assignments which maximize the number of conserved interactions between two species. For that purpose they use the following probabilistic model. They associate a binary variable $z_{ij}$ to each possible protein ortholog pair $(f_i,y_j)$ (here $f_i$ and $y_j$ denote Fly and Yeast proteins from the same Inparanoid cluster), where $z_{ij}=1$ means that $f_i$ and $y_j$ are functional orthologs. Two variables $z_{ij}$ and $z_{kt}$ are connected if at least one pair of proteins $(f_i,f_k)$ or $(y_j,y_t)$ is connected in its PPI network, and the other one has a common neighbor (or is also connected).Let $N(ij)$ denote the set of indices connected to $z_{ij}$. Then the probability law of $z_{ij}$  is modeled by:
\begin{equation}
P(z_{ij}|z_{N(ij)})=\frac{1}{1+exp\{-\alpha-\beta \sum_{kt \in N(ij)}z_{kt}\}}\,.
\end{equation}
The interpretation of this formula is that $z_{ij}$ has more chances to be equal to one when the number of neighbors equal to one is large. When there are only two proteins in cluster $f_i$ and $y_j$ then by definition $z_{ij}=1$. If $f_i$ and $y_j$ are from different clusters then also by definition $z_{ij}=0$. The parameters $\alpha$ and $\beta$ are estimated on the basis of training data, then a Gibbs sampling is performed to define the value of unknown variables z on the test set. We refer to \cite{Bandyopadhyay2006Systematic} for more details on this method.

\item \emph{MP method for exact optimization.} Although intractable in general, we now show that constrained GNA problem (\ref{eq:constrainedGNA}) can be solved exactly and efficiently in some cases, and propose a new, efficient algorithm based on message passing for that purpose. More precisely, we consider the situation where the set of proteins have been clustered into a finite set of $L$ groups $c_1,\ldots,c_L$, which form a partition of $V_G\cup V_H$, and where only proteins within the same group can be matched\footnote{Technically, we add dummy nodes in each cluster to obtain the same number of proteins of each species in each cluster.}. This situation, illustrated in Figure \ref{fig:inp_clust_schema}, represents for example the problem investigated by \cite{Bandyopadhyay2006Systematic}, where proteins of two organisms are first clustered by the Inparanoid algorithm, and functional orthologs are searched within clusters.
\begin{figure}[htbp]
 \begin{center}
\includegraphics[width=0.45\textwidth]{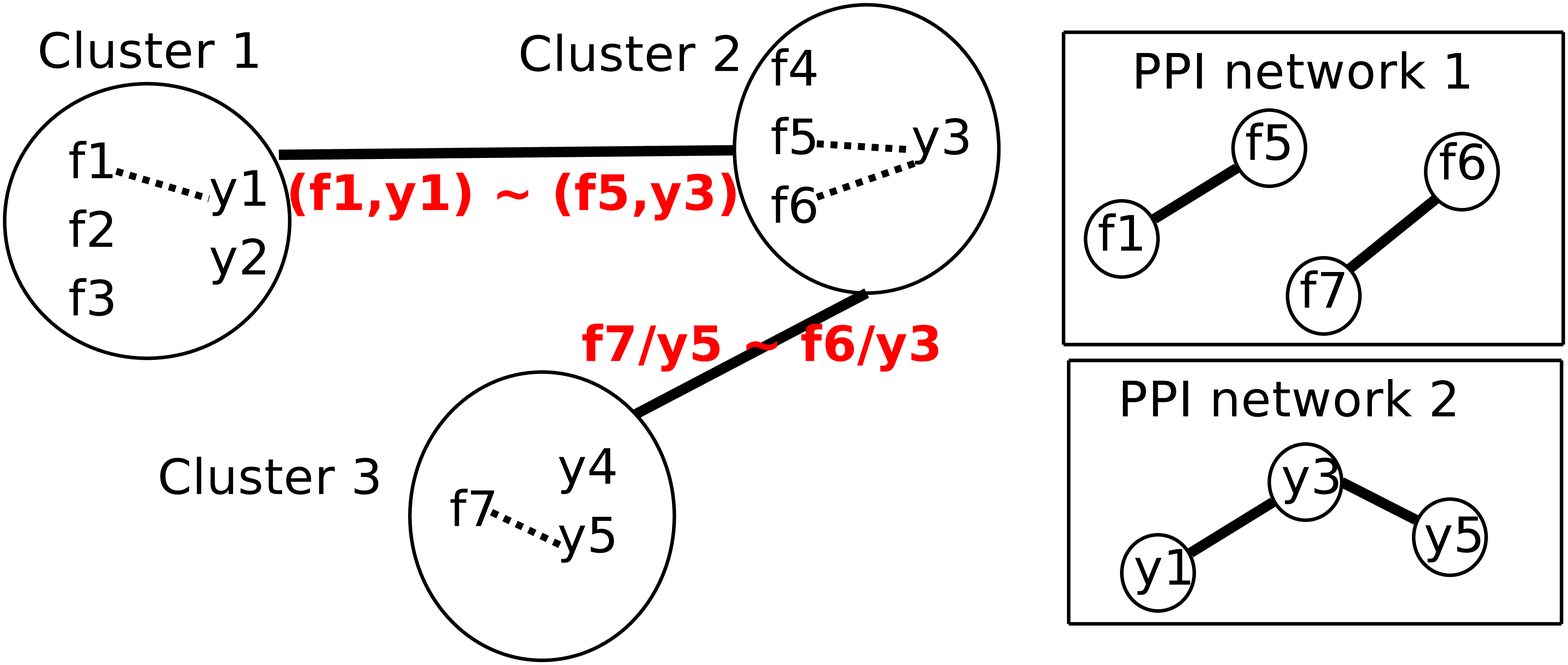}
\caption{Inparanoid cluster network. Two clusters are connected if there exist at least one pair of proteins in one cluster, and one pair of proteins in the other cluster, which may produce a conserved interaction.}
\label{fig:inp_clust_schema}
\end{center}
\end{figure}  
Let us now consider the $L$ clusters as vertices of a graph, and connect two clusters $c_i$ and $c_j$ if they contain proteins of both organisms that interact in their respective PPI network. For example, in Figure \ref{fig:inp_clust_schema}, $c_1$ and $c_2$ are connected because $c_1$ contains $f_1$ from the first organism and $y_1$ from the second organism, which interact with $f_5$ and $y_3$ respectively, both in $c_2$. 
The reason why we introduce this graph of clusters is that it allows to decompose the choice of a global matching $P$ into local matchings within each cluster, the dependency between the local choices being described by the edges of the graph. For example, if a cluster is isolated, then the choice of the matching within this cluster has no influence over the total number of conserved interactions apart from interactions within this cluster.  In other words, the local matching within an isolated cluster can be optimized independently from the others. On the other hand, if a cluster is connected to other clusters, then changing the matching within this cluster can affect the total number of interactions between proteins of different clusters, and the matchings between connected clusters must be chosen synchronously to optimize the total number of conserved interactions.

More formally, if we denote by $P_1,\ldots,P_L$ the permutation $P$ restricted to the $L$ clusters, then an important property is that the total number of interactions conserved by $P$ decomposes as:
\begin{equation}\label{eq:J_factor}
J(P) = \sum_{i=1}^L J_1(P_i) + \sum_{i\sim j} J_2(P_i,P_j)\,,
\end{equation}
where $J_1(P_i)$ denotes the number of conserved interactions within $c_i$, $J_2(P_i,P_j)$ denotes the number of conserved interactions between $c_i$ and $c_j$,  and $i\sim j$ means that $c_i$ is connected to $c_j$.

While maximizing (\ref{eq:J_factor}) remains a challenging optimization problem in general, it may be optimized efficiently if the graph of clusters has a particular structure, e.g., if many nodes are isolated or if it contains no loop. For example, Figure \ref{fig:inp_complex}(a) shows the graph of clusters for the problem of fly/yeast protein alignment investigated by \cite{Bandyopadhyay2006Systematic}. Interestingly, this graph has no loop.  In this case, we can maximize  (\ref{eq:J_factor}) by a particular Message Passing (MP) algorithm \cite{Jordan2001Learning}. The idea of the MP  algorithm is similar to the Viterbi algorithm \cite{Viterbi1973Error} widely used to optimize functions over linear graphs, such as finding the most likely set of hidden states in a hidden Markov model \cite{Durbin1998Biological}. Here we describe how to apply MP on a graph without loop to optimize (\ref{eq:J_factor}). First, we note that each of the permutations involving proteins within a connected component of the graph can be optimized independently from each other, so we just consider a single connected component without loop, i.e., a tree $\Tcal$ of clusters. We choose a vertex of $\Tcal$ that we call root, which allows to define the directions up (towards the root) or down (away from the root) when moving on edges of the graph. Each cluster $c_i$ except the root has a unique parent cluster, namely, the connected cluster in the direction of the root. The clusters connected to a cluster $c$ which are not its parent are called its children and are denoted $ch(c)$. To each node $c$ of $\Tcal$, we associate a vector $u_c \in \RR^{\Pcal_c}$, where $\Pcal_c$ is the set of possible local matchings within $c$, i.e., the set of possible $P_c$'s. The MP algorithm to solve (\ref{eq:J_factor}) is then a recursive algorithm, which starts from the leaves up to the root in a first phase (the ``forward'' step) to find the optimal value of the functional, and then downwards from the root to leaves (the ``backward'' step) to find the solution which achieves the optimal value. The forward step at node $c$ solves, for any $P_c \in \Pcal_c$:
\begin{equation}\label{eq:mp}
u_c(P_c) = J_1(P_c) + \sum_{c' \in ch(c)} \max_{P_{c'} \in \Pcal_{c'}} \left[ u_{c'}(P_{c'}) + J_2(P_c,P_{c'})\right]\,.
\end{equation}
At the end of the forward step, the maximum value of the vector $u$ at the root is equal to the maximal value of $J(P)$, and the local permutation which achieves this maximum is the optimal local permutation. In the backward step, the optimal local matching of the children of a cluster are obtained by recovering the local permutations $P_{c'}$ which achieved the optimal value in (\ref{eq:mp}) for the optimal permutation of the parent cluster.

We note that it is also possible to use the MP algorithm on graphs that are not trees, but which have a small tree-width value \cite{Jordan2001Learning}. Roughly speaking it means that the graph of clusters is not a tree, we may transform it into a tree by grouping together clusters. If the size of these cluster groups is not very large, then the exact optimization may still be feasible. 
\end{itemize}

\section{Data}
\label{sec:data}
In order to compare the performance of the different graph matching methods, we performed several experiments aiming at aligning the PPI networks of the yeast \emph{S. cerevisiae} and of the fly \emph{D. melanogaster}, as already investigated by \cite{Bandyopadhyay2006Systematic} and \cite{Singh2008Global}. We downloaded all necessary data from the supplementary materials of \cite{Bandyopadhyay2006Systematic}\footnote{{\tt http://www.cellcircuits.org/Bandyopadhyay2006}{http://www.cellcircuits.org/Bandyopadhyay2006}}. \OMIT{\todo{Describe the data: how many proteins, interactions, where do they come from, measure of similarity...}.} The yeast PPI network contains 4,389 proteins and 14,319 pairwise interactions, while the fly network contains 7,038 proteins and 20,720 interactions. In addition we also retrieved the set of Inparanoid clusters used by \cite{Bandyopadhyay2006Systematic}, consisting in 2,244 cluster covering 2,834 yeast proteins and 3,881 fly proteins. The majority of these clusters (1,552) contains only two proteins (one from fly, one from yeast), while the remaining 692 cluster contain at least two proteins from the same species and one from the other species. Those 692 clusters are called ambiguous in \cite{Bandyopadhyay2006Systematic}, since they do not allow to associate a single protein from the fly to a single protein from the yeast as functional orthologs.
\OMIT{ and we need  additional information to resolve these ambiguous association. The information on PPI networks may be used in 121 of these clusters as it was stated in \cite{Bandyopadhyay2006Systematic}.}

\section{Results}
\label{sec:results}
We wish to investigate two different questions: (i) compare the ability of the different methods to find alignment with many conserved interactions, and (ii) assess whether conserving more interactions really helps in retrieving more functional orthologs. While the first question can be answered without ambiguity by counting the number of conserved interactions found by the different methods in different settings, the second one, as we will see, remains difficult to answer due to the lack of large-scale and curated ground truth.

We performed three sets of experiments, in order to compare the different methods in different settings and to test different formulations of the GNA problem. In the first set of experiments, we reproduce the problem studied by \cite{Bandyopadhyay2006Systematic}, where the goal is to disambiguate functional orthologs within Inparanoid clusters using PPI information. This is a particular instance of the constrained GNA problem which turns out to be amenable to exact optimization by the MP method. In the second set of experiments, we generalize the benchmark problem of \cite{Bandyopadhyay2006Systematic} by adding second-order interactions between proteins in order to account for possible noise in the interaction data or protein duplications. In that case we are again confronted with a constrained GNA problem, but the increased number of interactions makes its exact minimization intractable and only approximate methods for constrained GNA can be applied. Finally, in a third set of experiments, we discard the knowledge of Inparanoid clusters and directly search a global alignment which balances the similarity between aligned proteins and the number of conserved interactions. This is then an instance of the balanced GNA problem. In all cases, we assess the number of conserved interactions captured by the different methods, as an indicator of how well they solve the GNA problem. Furthermore, since the final objective of PPI network alignment is to match functional orthologs, we assess for each method how many matched pairs are present in the HomoloGene database, a set of curated functional orthologous pairs based on the comparison of the protein as well as the DNA sequence which we consider here as a "gold standard" for disambiguation purpose.

\subsection{Disambiguation of functional orthologs within Inparanoid clusters}
\label{sec:constrained_first_order}
\begin{figure*}[htbp]
 \begin{center}
(a)\includegraphics[width=0.45\textwidth]{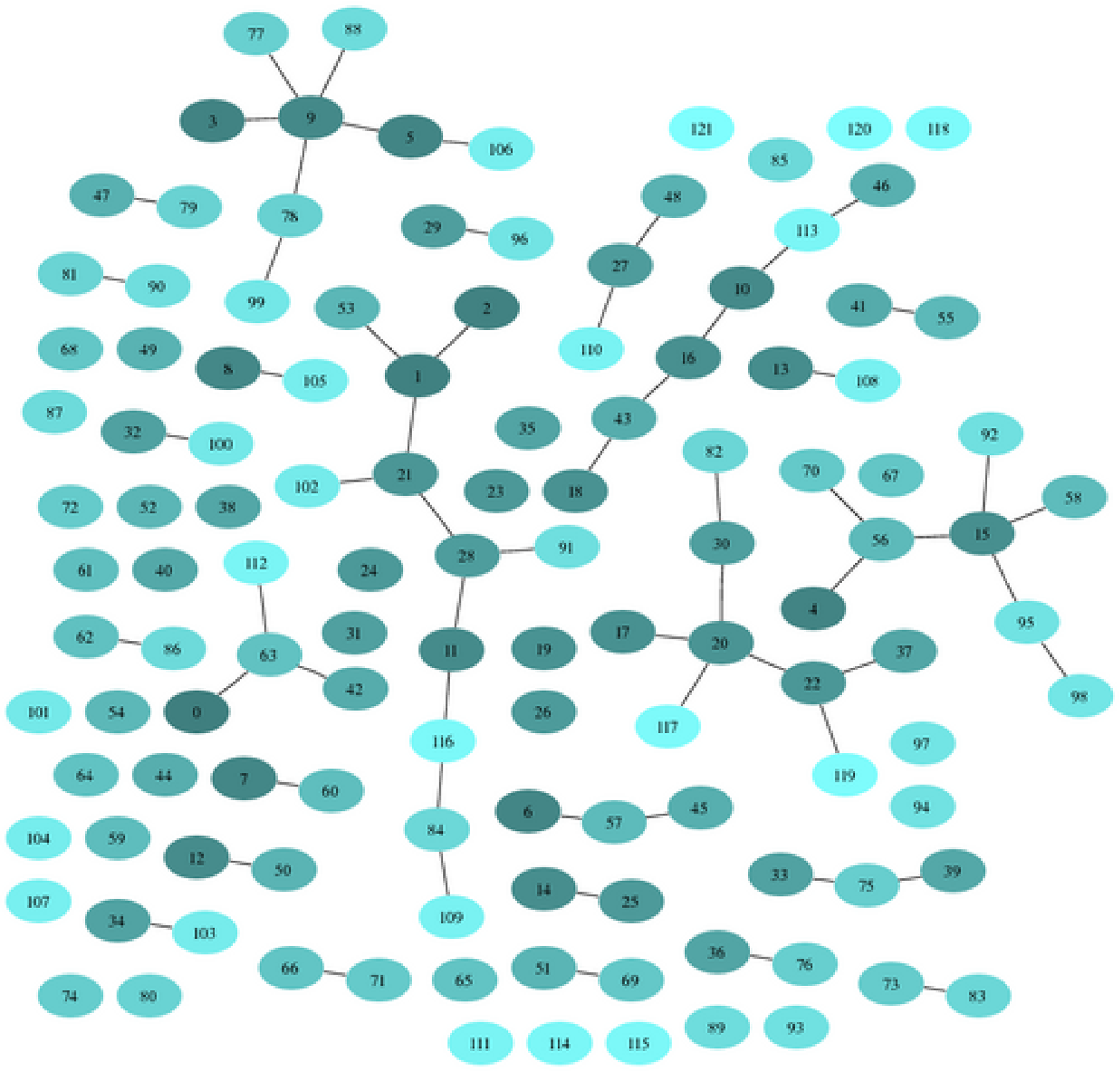}(b)\includegraphics[width=0.45\textwidth]{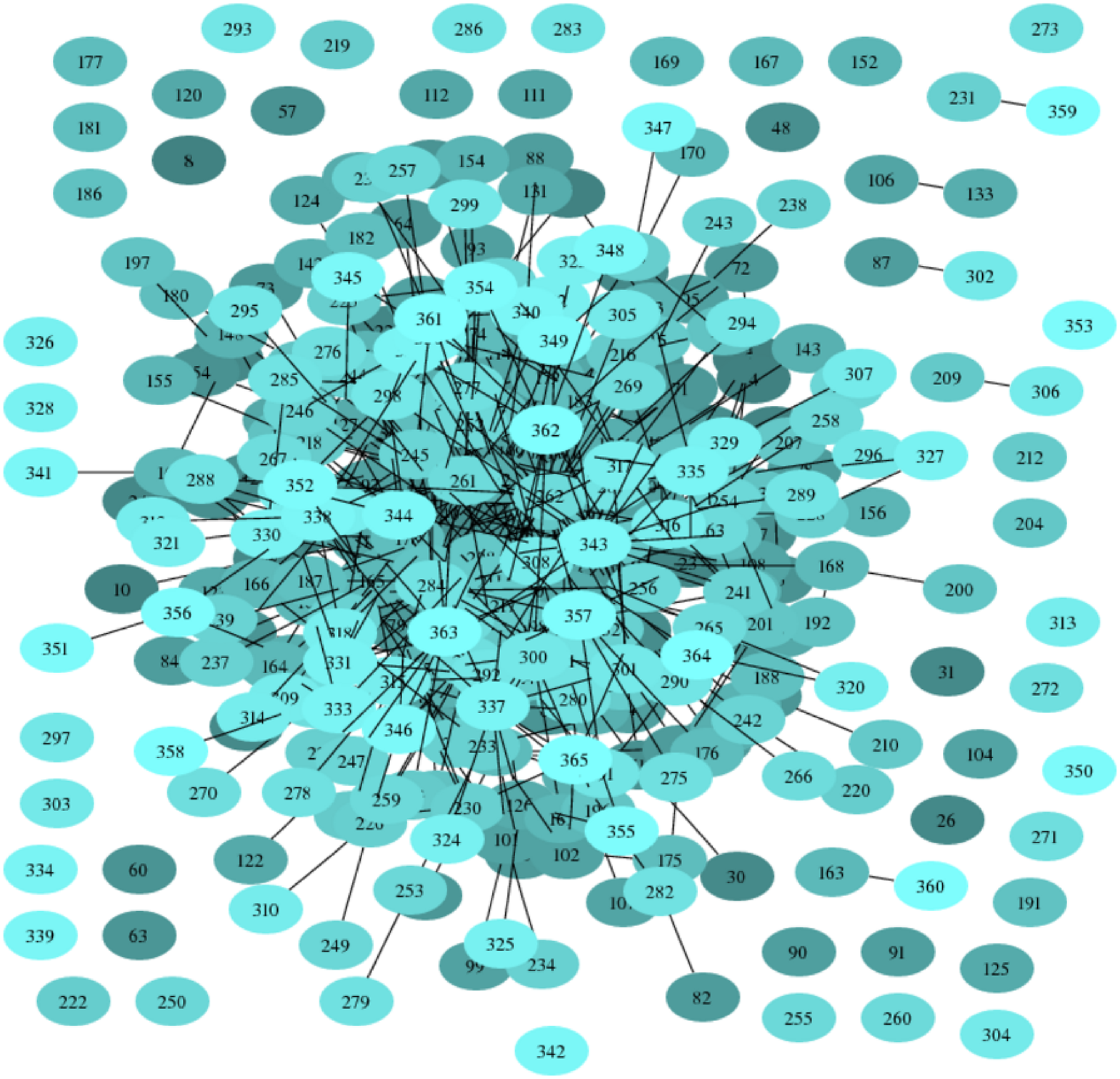}
\caption{Inparanoid cluster networks. (a) The case of the benchmark data used in \cite{Bandyopadhyay2006Systematic}. (b) The  case of generalized interactions (1-4), see text.}
\label{fig:inp_complex}
\end{center}
\end{figure*}

The goal of this experiment  is to use PPI GNA to select functional orthologs between the yeast and the fly for proteins with several homologs. More precisely, all proteins sequences are first clustered into groups by the Inparanoid algorithm \cite{Brein2005Inparanoid}, and only proteins from the same cluster can be considered as protein functional orthologs. Then each GNA algorithm tries to find an association of protein functional orthologs which maximizes the total number of conserved interactions. In other words, we try to solve the constrained GNA (\ref{eq:constrainedGNA}), where the constraints are provided by the Inparanoid clusters. A priori, the most natural definition of ``conserved interaction'' for the alignment $(f_1-y_1)$ and $(f_2-y_2)$ (where $f_1$ and $f_2$ are fly's proteins, and $y_1$ and $y_2$ are yeast's proteins) is the following:
\begin{enumerate}
\item $f_1$ interacts with $f_2$, and $y_1$ interacts with $y_2$ in their respective PPI networks.
\end{enumerate}
However, this strict notion of conserved interaction leads to a very small number of potentially conserved interactions. To have more potential interactions, \cite{Bandyopadhyay2006Systematic} generalized this definition by adding the following two cases, which additionally allow to account for possible duplication or fusion events in the two proteomes:
\begin{enumerate}
\setcounter{enumi}1
\item $f_1$ interacts with $f_2$ in the fly PPI network, and $y_1$ has a common neighbor with $y_2$ in the yeast PPI networks;
\item $f_1$ has a common neighbor with $f_2$ in the fly PPI network, and $y_1$ interacts with $y_2$ in the yeast PPI networks.
\end{enumerate}
To be able to compare the results of different algorithms, we use this exact definition of conserved interactions (cases 1-3). Figure \ref{fig:inp_complex}(a) presents the network of Inparanoid clusters (as explained in Figure \ref{fig:inp_clust_schema}) used in \cite{Bandyopadhyay2006Systematic}, where only non-isolated ambiguous clusters are shown. As can be easily seen, this network which contains 121 ambiguous clusters has no loop, which implies that we can use the MP method to find the optimal alignment with the largest number of conserved interactions.
Although we know how to solve the problem exactly in this case with the MP method, it is instructive to compare also the results of the different approximate algorithms for constrained GNA, namely, MRF and the constrained versions of IsoRank, GA and PATH. 
\OMIT{We also added a baseline method, denoted BLAST, which does not use the PPI network information to find the optimal alignment. The BLAST method aligns proteins of two species on the basis of their sequence similarities, by maximizing the sum of all similarity BLAST scores BLAST(i,j) over all pairs (i,j) of associated protein functional orthologs. }
To construct the alignment made by the MRF method \cite{Bandyopadhyay2006Systematic}, we downloaded the result file\footnote{http://www.cellcircuits.org/Bandyopadhyay2006/data/Bandyopadhyay\_results.xls} with probabilities for all possible protein association, and we extracted the one-to-one alignment by taking the most probable pairs. The results of the PATH, GA and IsoRank algorithms were obtained with the GraphM package \cite{Zaslavskiy2008GraphM}. 

Table \ref{tab:res} presents the results of all algorithms on this benchmark, in terms of conserved interactions, number of HomoloGene pairs, and running time.
 \begin{table*}
\centering
\caption{Performance of the different methods for constrained GNA on the benchmark of \cite{Bandyopadhyay2006Systematic}. \OMIT{\todo{Add more explanations here.}}Each algorithm is evaluated by the number of conserved interactions, number of recovered HomoloGene pairs and the running time. The number of recovered HomoloGene pairs is counted only in 121 ambiguous Inparanoid clusters where PPI data may be used.}
\begin{tabular}{|c|c|c|c|c|c|}
\hline
Algorithm & MP & MRF & IsoRank & GA & PATH \\
\hline
Number of conserved interactions & \bf 238 & 233 & 228 & \bf 238 & \bf 238 \\
\hline
Number of HomoloGene pairs (121 cl.) & 41 & 36 & 39 & 41&41\\
\hline
Timing(sec) & 1-2  & 10 &1-2&1-2&80-100\\
\hline
\end{tabular}
\label{tab:res}
\end{table*}
\OMIT{We first observe that, not surprisingly, all algorithms which use the information of the PPI network structure significantly increase the number of conserved interactions over the baseline BLAST method.}
We know that the MP algorithm produces the maximal possible value (238 in this case), and an interesting observation is that  the GA and the PATH algorithms reach this maximum, while the MRF (233) and the IsoRank (228) algorithms do not. All methods are comparable in terms of CPU time, except for MRF which is one order of magnitude slower on this dataset. Although the differences in number are slight, with only $2\%$ more conserved interactions for MP/GA/PATH than for MRF, and $4\%$ more than for IsoRank, this nevertheless confirms that even on this relatively easy optimization problem neither MRF nor IsoRank finds the optimal solution, which can be found by other methods at no additional computational cost.

\begin{figure*}[htbp]
 \begin{center}
(a)\includegraphics[width=0.3\textwidth]{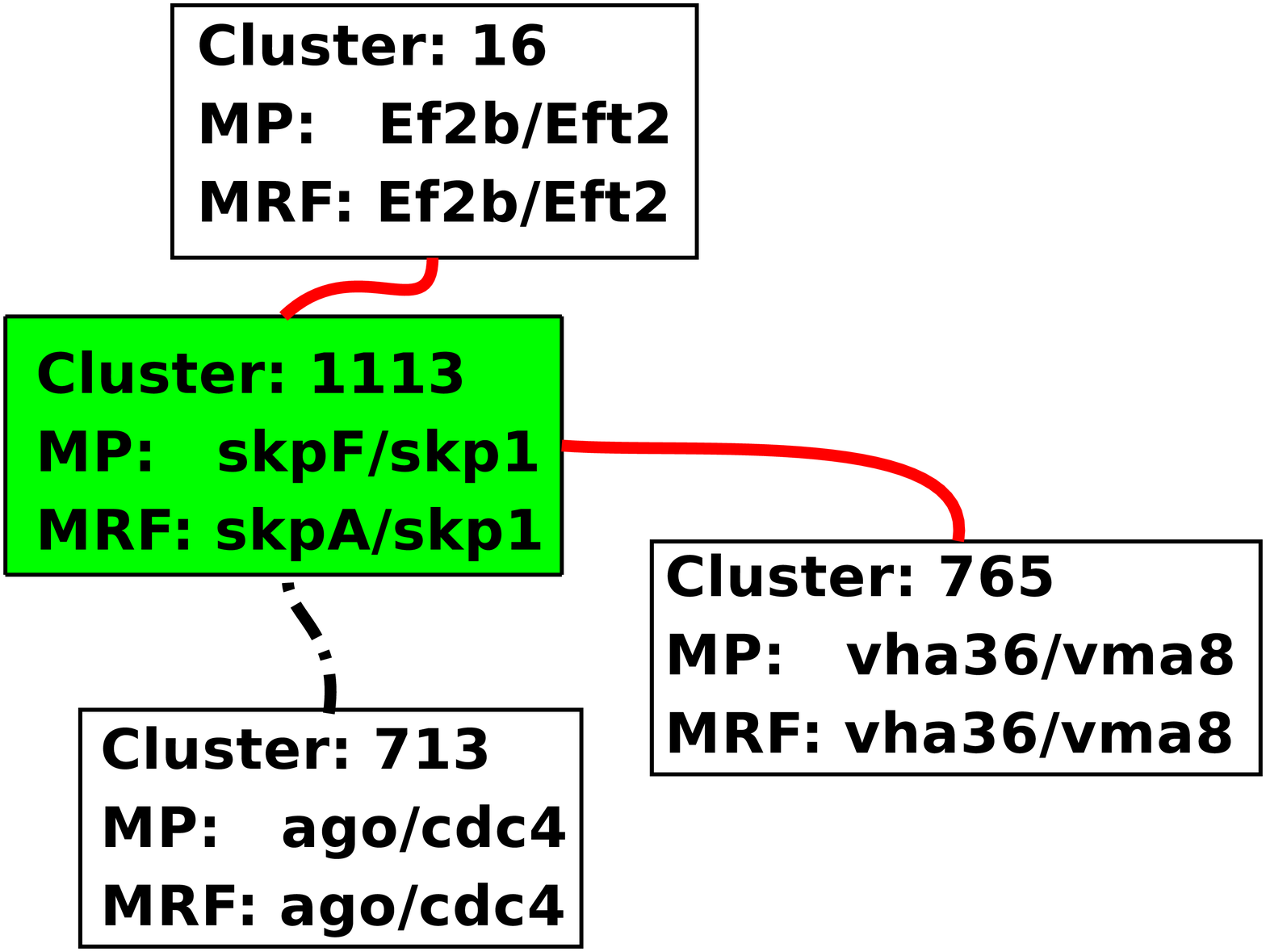}
\qquad (b)\includegraphics[width=0.55\textwidth]{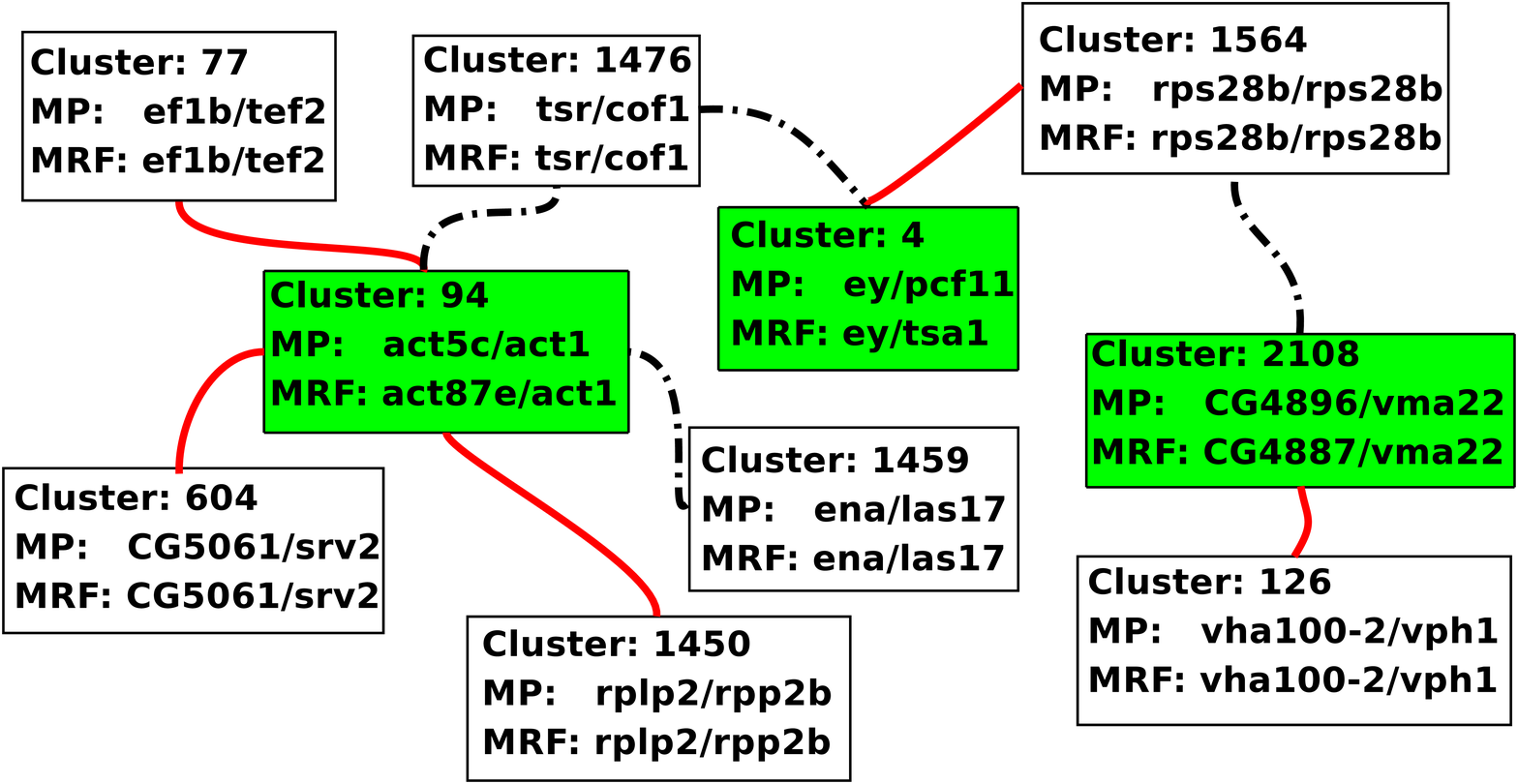}
\caption{Illustration of difference between MRF and MP alignment. Each box represents an Inparanoid cluster, white unfilled boxes represent clusters where MP and MRF assignments are the same. Red solid lines represent interactions conserved by MP alignment and not by MRF, black dotted lines represent interactions conserved by MRF and not by MP.}
\label{fig:diff}
\end{center}
\end{figure*}  

Figures \ref{fig:diff}(a) and \ref{fig:diff}(b) show some examples where the MRF assignment and the assignment made by the MP, PATH and GA algorithms are different, and illustrate how these differences influence the total number of conserved interactions. For instance, in the Inparanoid cluster 1113, the MRF algorithm associate the fly protein skpA to the yeast protein skp1, while the MP algorithm prefers the assignment skpF to skp1. In the later case we lose one conserved interaction with pair ago-cdc4, but we gain two new conserved interactions with (vha36,vm28) and (ef2b,eft2). In another example, shown in Figure \ref{fig:diff}(b), the MP algorithm proposes a different association for the yeast protein act1 in the 94-th Inparanoid cluster. This assignment results in two lost and three gained conserved interactions. From a biological point of view, the assignment of the fly protein act87e to act1 proposed by the MRF algorithm seems to be worse that the assignment (act5c,act1) proposed by the MP algorithm. Indeed, although proteins act5c and act87e are very similar (being both from the actine family), it is known that act1 and act5c participate together to the INO80 protein complex (which exhibits chromatin remodeling activity and 3' to 5' DNA helicase activity), while act87e does not.

In order to assess more systematically and quantitatively whether differences in the number of conserved interactions lead to significant differences in number of correctly assigned functional orthologous pairs, we counted how many pairs in each alignment is reported as functional orthologous in the HomoloGene database, considered here as a "gold standard". As shown in Table \ref{tab:res}, the number of HomoloGene pairs in each alignment also differs between the different methods, ranging from 36 for MRF to 39 for IsoRank and 41 for MP/GA/PATH. Interestingly, we observe that the method MP, GA and PATH, which retrieve the largest number of conserved interaction, also result in the largest number HomoloGene pairs (41), which represents a relative increase of $13\%$ compared to MRF (36), and of $5\%$ compared to IsoRank. To illustrate the differences between the methods, Table \ref{tab:homologene} lists the HomoloGene pairs found by MRF and not MP/GA/PATH, and vice versa. Interestingly, a new  method for PPI network alignment was published recently \cite{Yosef2008Improved}, which detects 37 HomoloGene orthologs on the same set of proteins. This puts its between MRF and IsoRank according to this criterion.
\begin{table}
\centering
\caption{HomoloGene orthologs found by the MP method and not by MRF and vice versa.}
\begin{tabular}{|c|c||c|}
\hline
\multicolumn{2}{|c||}{MP}&MRF\\
\hline
(TfIIA-S, TOA2)&(RPL23, RPL23A)& (Pros35, PRE5)\\
(CG13890, ECI1)&(Gapdh1, TDH1)& (Rab11, Ypt31)\\
(TfIIS, DST1)&(Rpt4, Rpt4)& (Rps26, Rps26A)\\
(Ef1gamma, TEF4)&(act5c, act1)&(CG6523, YDR098C)\\
(Glut1, YBR241C)&(Sir2, hst1)&(CG8690, YBR299W)\\
\hline
\end{tabular}
\label{tab:homologene}
\end{table}

The validity of taking HomoloGene as a "gold standard" for assessing the number of correctly assigned homologous pairs remains, however, subject to discussion. Indeed, although HomoloGene clusters are defined using a variety of evidences, they are mainly driven by sequence similarity. To illustrate this, we assessed the performance of a simple alignment method which matches pairs within an ambiguous cluster by maximizing the total sequence similarity over matched pairs. This method does not use any PPI information for the matching. The resulting alignment has only 184 conserved interaction, which is not surprisingly much worse than all methods which take PPI into account. However, the resulting matched pairs contain 43 HomoloGene pairs, which is more than all methods taking into account PPI. This shows that the number of HomoloGene pairs as an indicator should be taken with caution, since it favors methods which focus on matching proteins based on sequence similarity only.

\subsection{Disambiguation of Inparanoid clusters with second-order interactions} 
\label{sec:constrained_sec_order}
The idea of \cite{Bandyopadhyay2006Systematic} to expand the natural notion of conserved interaction (case 1) to cases 2 and 3, aims to take into account second-order interactions, that is, when two proteins do not interact directly to each other have a common neighbor. Another natural generalization of the notion of conserved interaction is then the following case:
\begin{enumerate}
\setcounter{enumi}3 
\item $f_1$ has a common neighbor with $f_2$, and $y_1$ has a common neighbor with $y_2$, in their respective PPI networks.
\end{enumerate}
Adding interactions according to this rule makes the problem computationally more difficult, since ambiguous clusters become more connected. Indeed, while we were able to solve the original problem exactly with the MP algorithm, the network of Inparanoid clusters when cases 1-4 are included takes the form presented in Figure \ref{fig:inp_complex}(b). Contrary to the previous network (cases 1-3 in Figure \ref{fig:inp_complex}(a)), the new network has loops and is not amenable to exact optimization with the MP procedure. Only approximate algorithms can be applied in this case.

In order to compare all methods (except MP) in this new setting, we re-implemented the MRF algorithm with the new data. The estimated values of the model parameters (see details in \cite{Bandyopadhyay2006Systematic}) are $(\alpha=0.51,\beta=-6.87)$. We used the same training and test data as those used used in \cite{Bandyopadhyay2006Systematic} to estimate them. Then we estimated the probabilities of being protein orthologs for potential pairs of proteins by Gibbs sampling, and obtained a one-to-one alignment based on the most probable associations.

Table \ref{tab:res2} shows the results obtained by the different graph matching algorithms. Although we do not know the maximum number of interactions that can be conserved in this case, we observe again that PATH and GA find solutions with $3-4\%$ more interactions conserved than MRF and IsoRank. There is no clear difference in the number of HomoloGene pairs between the different methods, and the addition of second-order interactions has no obvious effects on this indicator neither: it leads to a gain of 3 pairs for MRF, but to a loss of one pair for IsoRank and PATH, and to no change for GA.
\OMIT{
\begin{table*}
\centering
\caption{Performance of the different methods for constrained GNA on the benchmark of \cite{Bandyopadhyay2006Systematic} with second-order interactions added. The number of recovered HomoloGene pairs is counting on the 121 Inparanoid clusters from the previous section as well as on the new 602 ambiguous Inparanoid clusters have second-order interaction with other Inparanoid clusters}
\begin{tabular}{|c|c|c|c|c|c|c|}
\hline
Algorithm & MRF & IsoRank & GA & PATH & BLAST\\
\hline
Number of conserved interactions & 1112&1101&1140&1143&804\\
\hline
Number of HomoloGene pairs (121 cl.) &39 & 38&41 &40 &43\\
\hline
Number of HomoloGene pairs (602 cl.) &172 & 167&172 &166 &191\\
\hline
Timing(sec) & 623 &31&372&1542&2-3\\
\hline
\end{tabular}
\label{tab:res2}
\end{table*}
}

\begin{table*}
\centering
\caption{Performance of the different methods for constrained GNA on the benchmark of \cite{Bandyopadhyay2006Systematic} with second-order interactions added. The number of recovered HomoloGene pairs is counting on the 121 Inparanoid clusters from the previous section as well as on the new 602 ambiguous Inparanoid clusters have second-order interaction with other Inparanoid clusters}
\begin{tabular}{|c|c|c|c|c|}
\hline
Algorithm & MRF & IsoRank & GA & PATH \\
\hline
Number of conserved interactions & 1,112&1,101&1,140&1,143\\
\hline
Number of HomoloGene pairs (121 cl.) &39 & 38&41 &40 \\
\hline
Number of HomoloGene pairs (602 cl.) &172 & 167&172 &166 \\
\hline
Timing(sec) & 623 &31&372&1,542\\
\hline
\end{tabular}
\label{tab:res2}
\end{table*}

\subsection{Global PPI network alignment by balancing sequence and interaction conservation}
In this last series of experiments, we consider the problem proposed by \cite{Singh2008Global}, for which IsoRank reflects the state-of-the-art: find a global PPI alignment by balancing the sequence similarity in matched pairs with the total number of conserved interactions, allowing in particular matches between proteins in different Inparanoid clusters if they allow an increased number of conserved interactions. For this application we can only compare the three methods for balanced GNA, namely, IsoRank, GA and PATH. The trade-off between matching proteins with similar sequences and matching with a lot of conserved interactions is controlled by the parameter $\lambda$ in (\ref{eq:balancedGNA}) and (\ref{eq:R2}). The greater $\lambda$, the more attention we pay to the sequence similarity and the less to the number of conserved interactions. For each method, by varying $\lambda$, we therefore obtain a family of alignments with different compromise found between the number of conserved interactions $J(P)$ (\ref{eq:balancedGNA}) and the summary sequence similarity score $S(P)$ (\ref{eq:balancedGNA}).  

\begin{figure}[htbp]
 \begin{center}
\includegraphics[width=0.48\textwidth]{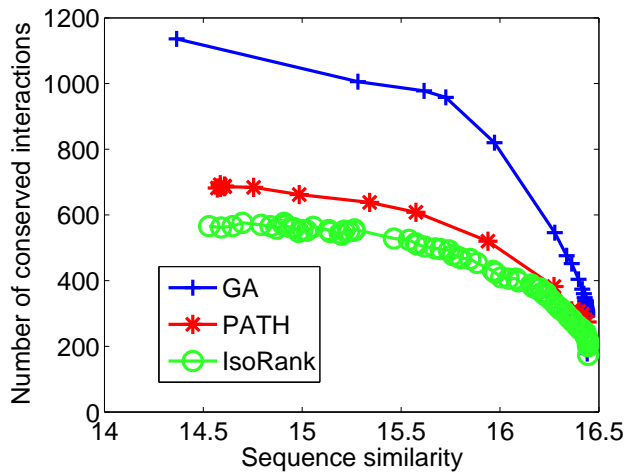}
\caption{Algorithm performance comparison. Number of conserved interaction $J(P)$ versus sequence similarity $S(P)$.}
\label{fig:num_cons_inter_noconstraints}
\end{center}
\end{figure} 
Figure \ref{fig:num_cons_inter_noconstraints} shows the different trade-offs which are found by the different methods. For a given level of average sequence similarity, we wish to have the largest possible number of conserved pairs. We observe that over all the range of average sequence similarity, the GA algorithms clearly outperforms PATH, which itself outperforms IsoRank. For example, for the trade-off parameter choice advocated by \cite{Singh2008Global} for IsoRank ($\lambda=0.6$), IsoRank finds an alignment with 566 conserved interactions, corresponding to an average sequence similarity score in the matched pairs of $15.26$. At this level of average sequence similarity, PATH and GA find alignments with respectively $678$ and $1,006$ interactions, which corresponds to relative improvements of respectively  $20\%$ and $78\%$.

Again, there is still only limited objective evidence that optimizing the number of conserved interactions leads to better matching in terms of functional orthology detection. As an attempt to test this fact, we first counted, for each alignment, the number of HomoloGene pairs in the alignment. However, we observed that, for each method, this number increases monotonically when more weight is given to sequence similarity as opposed to interaction conservation. This again highlights the limitation of this criterion, which is optimized by construction when sequences are optimally matched in terms of similarity. We then attempted to compare the different alignments in terms of mean similarity between gene ontology (GO) annotations of matched pairs. In order to compare GO annotations of two proteins we tested the method presented by \cite{Singh2008Global} to compute the functional coherence of a pair. However, we were not able to observe any clear difference between the methods, or between the different parameter choice for each individual method. The maximum mean functional coherence over the choice of the trade-off parameter is respectively $0.519$, $0.509$ and $0.522$ for IsoRank, GA and PATH. However the fluctuations of this score when the parameters change are so large that these maximum values are not significantly different. This is due to the fact that the number of annotated proteins remains limited, and that they are rarely annotated with such precision that it is possible to clearly differentiate true functional orthologs from spurious ones \cite{Bandyopadhyay2006Systematic}. For example, when we estimate the functional score of a given alignment, there is rarely more than $15-20\%$ of pairs with GO annotations.

\section{Discussion}
We presented two general formulations for the GNA problem. The constrained GNA formulation corresponds to a situation where we have a strong \emph{a priori} about which pairs can be matched. In the balanced GNA problem, we replace the binary constraints on which pairs are allowed by a more global objective function which balances the matching of similar proteins with the conservation of interactions, with a parameter to smoothly control the trade-off between these two contradictory goals. While MRF and IsoRank are popular methods for these two formulations, we proposed in this paper new methods which lead to significantly better alignments, when we assess the quality of an alignment in terms how many conserved interactions are retrieved. In particular, the MP method, when it is applicable, finds the optimal solution of a constrained GNA problem, and the GA method provides consistently good results in both cases. The question of which formulation is the best for a given application and dataset, between the constrained and balanced GNA, remains largely open and worth further systematic investigations. Regarding the relative performance of the different methods in terms of how many conserved interactions they find, we observed that the MP/GA/PATH methods outperform MRF and IsoRank in both situations. This is not so surprising given that, once the problem is explicitly stated as a graph matching problems, it makes sense to use methods borrowing ideas and techniques from state-of-the-art graph matching approaches. The impressive performance of GA compared to PATH in the balanced GNA experiment (Figure \ref{fig:num_cons_inter_noconstraints}) is more surprising, given the good performance of PATH on a number of other benchmarks \cite{Zaslavskiy2008patha}. We believe this weakness of PATH is due to the large difference in the number of nodes between the two networks. Indeed, the resulting large number of dummy nodes that must be added generate singularities in the convex relaxation in the PATH algorithm. \OMIT{On the contrary, in the case of constrained GNA, the PATH algorithm is more efficient because most of the singularities disappears when we project the convex relaxation on the subspace defined by constraints.}  

The GNA problems we studied have several extensions. First, it may be interesting to consider alignment of weighted PPI networks with weights representing, for instance, experimental evidence of interaction existence. Interestingly, the  PATH, GA and IsoRank algorithm can be applied directly to a weighted network, by just replacing the binary graph adjacency matrix by a real-valued matrix. Another relevant extension is the alignment of multiple PPI networks, corresponding to more than two species, via pairwise comparisons as it was presented by \cite{Singh2008Global}. Finally, it may be relevant in some cases to match one protein of one species with several proteins of the other species, to account for possible duplications or fusion events. An interesting property of the PATH algorithm is the fact that estimate a permutation matrix by first solving a relaxed problem. The solution of the relaxed problem is a doubly stochastic matrix whose entries can be interpreted as probabilities for proteins to be functional orthologs \cite{Zaslavskiy2008patha} . Therefore, in order to allow many-to-many assignments of proteins, we could use the solution of the convex relaxation. 

Finally, although progresses in graph alignment algorithms can be monitored by objective quantitative measures such as the number of conserved interactions, their biological relevance remains difficult to assess. In particular, for the detection of functional orthologs, it is apparent that current GO annotations or curated databases of functional orthologs are either biased by construction (e.g., HomoloGene), or not precise enough and too scarce for systematic evaluation (e.g., GO annotations). We believe we are reaching a point where more experimental validations are needed. On the other hand, there are many other possible applications for efficient graph matching algorithms scaling to large biological networks, such as phylogenetic comparison of sets of networks, detection of new conserved pathways, or curation of PPI data. We expect the methods proposed in this paper to have a direct impact in these applications.

\bibliographystyle{plain}
%

\begin{thebibliography}{10}

\bibitem{Aebersold2003Mass}
R.~Aebersold and M.~Mann.
\newblock Mass spectrometry-based proteomics.
\newblock {\em Nature}, 422(6928):198--207, Mar 2003.

\bibitem{Almohamad1993linear}
H.A. Almohamad and S.O. Duffuaa.
\newblock A linear programming approach for the weighted graph matching
  problem.
\newblock {\em IEEE Trans. Inform. Theory}, 15(5):522--525, May 1993.

\bibitem{Bandyopadhyay2006Systematic}
S.~Bandyopadhyay, R.~Sharan, and T.~Ideker.
\newblock Systematic identification of functional orthologs based on protein
  network comparison.
\newblock {\em Genome Res.}, 16(3):428--435, Mar 2006.

\bibitem{Berg2006Cross-species}
J.~Berg and M.~L{\"a}ssig.
\newblock Cross-species analysis of biological networks by bayesian alignment.
\newblock {\em Proc. Natl. Acad. Sci. USA}, 103(29):10967--10972, Jul 2006.

\bibitem{Brein2005Inparanoid}
K.~Brein, M.~Remm, and E.~Sonnhammer.
\newblock Inparanoid: a comprehensive database of eukaryothic orthologs.
\newblock {\em Nucleic acids research}, 33, 2005.

\bibitem{Caelli2004eigenspace}
T.~Caelli and S.~Kosinov.
\newblock An eigenspace projection clustering method for inexact graph
  matching.
\newblock {\em IEEE Trans. Pattern Anal. Mach. Intell.}, 26(4):515--519, April
  2004.

\bibitem{Conte2004Thirty}
D.~Conte, P.~Foggia, C.~Sansone, and M.~Vento.
\newblock Thirty years of graph matching in pattern recognition.
\newblock {\em International journal of pattern recognition and artificial
  intelligence}, 18(3):265--298, 2004.

\bibitem{Durbin1998Biological}
R.~Durbin, S.~Eddy, A.~Krogh, and G.~Mitchison.
\newblock {\em Biological {S}equence {A}nalysis: {P}robabilistic {M}odels of
  {P}roteins and {N}ucleic {A}cids}.
\newblock Cambridge University Press, 1998.

\bibitem{Fields1989novel}
S.~Fields and O.~Song.
\newblock A novel genetic system to detect protein-protein interactions.
\newblock {\em Nature}, 340(6230):245--246, Jul 1989.

\bibitem{Flannick2006Graemlin}
J.~Flannick, A.~Novak, B.S. Srinivasan, H.H. McAdams, and S.~Batzoglou.
\newblock Graemlin: general and robust alignment of multiple large interaction
  networks.
\newblock {\em Genome Res.}, 16(9):1169--1181, Sep 2006.

\bibitem{Gold1996graduated}
S.~Gold and A.~Rangarajan.
\newblock A graduated assignment algorithm for graph matching.
\newblock {\em IEEE Trans. Pattern Anal. Mach. Intell.}, 18(4):377--388, April
  1996.

\bibitem{Jordan2001Learning}
Michael Jordan, editor.
\newblock {\em Learning in Graphical Models}.
\newblock The MIT Press, 2001.

\bibitem{Kelley2003Conserved}
B.P. Kelley, R.~Sharan, R.M. Karp, T.~Sittler, D.E. Root, B.R. Stockwell, and
  T.~Ideker.
\newblock Conserved pathways within bacteria and yeast as revealed by global
  protein network alignment.
\newblock {\em Proc. Natl. Acad. Sci. USA}, 100(20):11394--11399, Sep 2003.

\bibitem{Kelley2004PathBLAST}
B.P. Kelley, B.~Yuan, F.~Lewitter, R.~Sharan, B.R. Stockwell, and T.~Ideker.
\newblock {PathBLAST}: a tool for alignment of protein interaction networks.
\newblock {\em Nucleic Acids Res.}, 32(Web Server issue):W83--W88, Jul 2004.

\bibitem{Koyutuerk2006Pairwise}
M.~Koyut{\"u}rk, Y.~Kim, U.~Topkara, S.~Subramaniam, W.~Szpankowski, and
  A.~Grama.
\newblock Pairwise alignment of protein interaction networks.
\newblock {\em J. Comput. Biol.}, 13(2):182--199, Mar 2006.

\bibitem{Kuhn1955Hungarian}
H.~W. Kuhn.
\newblock The {H}ungarian method for the assignment problem.
\newblock {\em Naval Research}, 2:83--97, 1955.

\bibitem{Remm2001Automatic}
M.~Remm, C.E. Storm, and E.L. Sonnhammer.
\newblock Automatic clustering of orthologs and in-paralogs from pairwise
  species comparisons.
\newblock {\em J. Mol. Biol.}, 314(5):1041--1052, Dec 2001.

\bibitem{Sharan2005Conserved}
R.~Sharan, S.~Suthram, R.M. Kelley, T.~Kuhn, S.~McCuine, P.~Uetz, T.~Sittler,
  R.M. Karp, and T.~Ideker.
\newblock Conserved patterns of protein interaction in multiple species.
\newblock {\em Proc. Natl. Acad. Sci. USA}, 102(6):1974--1979, Feb 2005.

\bibitem{Singh2008Global}
R.~Singh, J.~Xu, and B.~Berger.
\newblock Global alignment of multiple protein interaction networks with
  application to functional orthology detection.
\newblock {\em Proc. Natl. Acad. Sci. USA}, 105(35):12763--12768, Sep 2008.

\bibitem{Sjoelander2004Phylogenomic}
K.~Sj{\"o}lander.
\newblock Phylogenomic inference of protein molecular function: advances and
  challenges.
\newblock {\em Bioinformatics}, 20(2):170--179, Jan 2004.

\bibitem{Suthram2005Plasmodium}
S.~Suthram, T.~Sittler, and T.~Ideker.
\newblock The plasmodium protein network diverges from those of other
  eukaryotes.
\newblock {\em Nature}, 438(7064):108--112, Nov 2005.

\bibitem{Umeyama1988eigendecomposition}
S.~Umeyama.
\newblock An eigendecomposition approach to weighted graph matching problems.
\newblock {\em IEEE Trans. Pattern Anal. Mach. Intell.}, 10(5):695--703, Sept.
  1988.

\bibitem{Viterbi1973Error}
A.~Viterbi.
\newblock Error bounds for convolutional codes and an asymptotically optimum
  decoding algorithm.
\newblock {\em IEEE Trans. Inform. Theory}, 13(2):260--269, 1973.

\bibitem{Yosef2008Improved}
N.~Yosef, R.~Sharan, and W.S. Noble.
\newblock Improved network-based identification of protein orthologs.
\newblock {\em Bioinformatics}, 24(16):i200--i206, Aug 2008.

\bibitem{Zaslavskiy2008GraphM}
M.~Zaslavskiy, F.~Bach, and J.~P. Vert.
\newblock {GRAPHM}: Graph matching package, 2008.
\newblock Available at \texttt{http://cbio.ensmp.fr/graphm}.

\bibitem{Zaslavskiy2008path}
M.~Zaslavskiy, F.~Bach, and J.~P. Vert.
\newblock A path following algorithm for graph matching.
\newblock In A.~Elmoataz, O.~Lezoray, F.~Nouboud, and D.~Mammass, editors, {\em
  Image and Signal Processing, Proceedings of the 3rd International Conference,
  ICISP 2008}, volume 5099 of {\em LNCS}, pp. 329--337. Springer Berlin /
  Heidelberg, 2008.

\bibitem{Zaslavskiy2008patha}
M.~Zaslavskiy, F.~Bach, and J.-P. Vert.
\newblock A path following algorithm for the graph matching problem.
\newblock Technical Report 00232851, HAL, 2008.
\newblock To appear in IEEE Trans. Pattern Anal. Mach. Intell.

\end{thebibliography}

\end{document}